\documentclass[10pt,a4paper]{article}

\usepackage{pifont,latexsym,amssymb,bbm}
\usepackage{yfonts,url}

\usepackage{hyperref}

\usepackage{graphicx}

\def\bd{\begin{description}}
\def\ed{\end{description}}

\def\beq{\begin{equation}}
\def\eeq{\end{equation}}
\def\bea{\begin{eqnarray}}
\def\eea{\end{eqnarray}}
\def\beas{\begin{eqnarray*}}
\def\eeas{\end{eqnarray*}}

\def\G1{\hbox{$\displaystyle{\mbox{\ding{172}}}$}}

\begin{document}

\title{Some paradoxes of  Infinity revisited}

\author{  Yaroslav D. Sergeyev\\
University of Calabria,  Rende, Italy\\    Lobachevsky University,
    Nizhni Novgorod, Russia\\   Institute of High Performance Computing and
  Networking \\of the National Research Council of Italy, Rende, Italy 
   }
 
 \date{}

\maketitle


\vspace{-8mm}

\begin{abstract}
In this article,  some classical paradoxes of infinity such as Galileo's paradox, Hilbert's paradox of the Grand Hotel, Thomson's lamp paradox, and the rectangle paradox of Torricelli are considered. In addition, three   paradoxes regarding divergent series and a new paradox dealing with multiplication of elements of an infinite set are also described. It is shown that the surprising  counting system  of an Amazonian tribe,  Pirah\~{a}, working with only three numerals (one, two, many) can help us to change our perception of these paradoxes. A recently introduced   methodology allowing one to work  with finite, infinite, and infinitesimal numbers in a unique computational framework  not only theoretically but also numerically is briefly described. This methodology is actively used nowadays   in numerous applications in pure and applied mathematics and computer science as well as in teaching. It is shown in the article that this methodology also allows one to consider the paradoxes listed above in a new constructive light.
\end{abstract}

\textbf{MSC:} 00A30, 97F30, 40-08, 40A05, 03A05, 97E40, 97C30.

\noindent  \textbf{Keywords: paradoxes of infinity, counting systems, Pirah\~{a}, grossone.}

\section{Introduction}

We use finite numbers every day and rarely think about \textit{the nature} of the infinite using it mechanically in our math classes. However, infinity and infinitesimals are among the most fundamental notions in mathematics (and not only). They
have attracted the attention of the most brilliant thinkers
throughout the whole history of humanity. Arabic, Indian, and
Babylonian mathematicians worked hard on these problems.
Aristotle, Archimedes, Euclid, Eudoxus, Parmenides, Plato,
Pythagoras, and Zeno  dealt with these problems in antiquity. In the years 1500--1900 important contributions were made
by such eminent research\-ers as Bolzano, Briggs, Cantor, Cauchy,
Dedekind, Descartes, Dirichlet, Euler, Hermite, Leibniz,
Lindemann, Liouville, Napier, Newton, Mercator, Peano, Stevin,
Wallis, and Weierstrass. In the $20^{th}$ century new exciting
results have been obtained by Brouwer, Cohen, Frege,
G\"{o}del, Hilbert, Robinson, Scott, and Solovay.

Introduction of the ideas of the number line, positional number
systems, negative numbers, zero, rational and irrational numbers,
limits, cardinal and ordinal numbers, continuum hypothesis, problems of consistency and
completeness, non-standard analysis are among the major milestones of these impressive
research efforts.  Research on these topics continues to be very active nowadays, as well (see, e.g., \cite{Alexander,Caldarola_paradoxes,Heller,Hellman_Shapiro_2013,Kanamori_book,Linnebo,Lolli_fil,Mancosu_2009,Mancosu_2016,
MM_bijection,Rizza_2,Ternullo_Fano} and references given therein).

However, it is well known that the ideas of infinities and infinitesimals lead to numerous paradoxes. Is it true that they are inevitable? Is it possible to propose a viewpoint allowing us to avoid some of them? In this paper, we try to answer these questions using counting systems of two tribes, Pirah\~{a} and Munduruk\'u, living in Amazonia nowadays (see \cite{Gordon,Pica}) together with a recent methodology working with finite, infinite, and infinitesimal
numbers in a unique computational framework not only theoretically but also numerically on a patented supercomputer called the Infinity Computer  (see a comprehensive technical survey \cite{EMS}, a brief survey in Italian \cite{UMI}, and a popular book \cite{Sergeyev} for its description).

\section{Paradoxes of infinity} Let us consider several    classical paradoxes  coming from different situations involving infinity. In many of them   the set, $\mathbb{N}$, of natural numbers
\beq
\mathbb{N} = \{1, 2, 3, 4,  5, \,\, \ldots \,\, \}
\label{paradoxes_N}
 \eeq
is involved. We informally define it as the set of numbers used to count objects.  Notice that nowadays not only positive integers are taken as elements of $\mathbb{N}$, but also zero is frequently included in $\mathbb{N}$. However, since historically zero has been invented significantly later with respect to positive integers used for counting objects, zero is not included in
$\mathbb{N}$ in this article.

\subsection{Galileo's paradox}  In his book ``Discourses and mathematical demonstrations relating to two new sciences'' published in 1638 Galileo Galilei  considered    the set $\mathbb{N}$ together with the set, that we call  $I^2$, of square natural numbers
 \beq
 I^2 = \{ x : x \in \mathbb{N}, i \in \mathbb{N}, x =i^{\,2} \}
 = \{1, 4, 9, 16, 25, \,\, \ldots \,\, \}.
 \label{4.3.0.1}
 \eeq
He then established the following bijection among the sets $I^2$ and the set of natural numbers,  $\mathbb{N}$, as
follows
 \beq
\begin{array}{ccccccc}
  \hspace{5mm} 1, & 2^2, & 3^2, & 4^2,  & 5^2, & 6^2, & \ldots    \\

  \hspace{5mm} \updownarrow &  \updownarrow & \updownarrow  & \updownarrow  & \updownarrow  &  \updownarrow &   \\

   \hspace{5mm}1, &  2, & 3, & 4 & 5,
       & 6,  &    \ldots \\
     \end{array}
\label{trad_square}
 \eeq
 This bijection is paradoxical since there are  much more numbers than squares and still to any number there can be found  the corresponding square and vice versa. Clearly, the same paradoxical result arises from considering a simpler bijection between $\mathbb{N}$ and the set  $\mathbb{E}$  of even numbers  being  a proper subset of $\mathbb{N}$:
 \beq
\begin{array}{ccccccc}
    \hspace{5mm} 2, & 4, & 6, & 8,  & 10, & 12, & \ldots    \\

  \hspace{5mm} \updownarrow &  \updownarrow & \updownarrow  & \updownarrow  & \updownarrow  &  \updownarrow &   \\

  \hspace{5mm}1, &  2, & 3, & 4 & 5,
       & 6,  &    \ldots \\
     \end{array}
\label{even_bij}
 \eeq
 From the modern point of view these bijections mean that all the sets involved are countable. However, the perplexity noticed by Galileo remains because in our every day life dealing with  finite objects and sets a part of a set is always less than the whole set.  In his \textit{Elements},
Euclid has expressed this property as Common Notion no.~5 `The whole is greater than the part', where Common Notions are evident assertions that are accepted without any proof.

Let us introduce now a new paradox that can be considered as a kind of inversion of (\ref{even_bij}) where we have started from $\mathbb{E}$ and established the bijection with the set  $\mathbb{N}$. In the new paradox, that hereinafter will be called \textit{set-multiplication paradox} we start from $\mathbb{N}$ and  arrive to  $\mathbb{E}$.

\subsection{Set-multiplication paradox}

Let us consider a finite even number $n$  and the corresponding set of natural numbers
 $$B = \{ 1, 2, 3, \ldots , n-2, n-1, n \}.$$
 Then we multiply each of its elements by 2 and obtain the set
   $$\bar{B} = \{ 2, 4, 6, \ldots ,  n-4, n-2, n, n+2, n+ 4, \ldots , 2n-4, 2n-2, 2n \}.$$
Notice the following three properties of the sets $B$ and $\bar{B}$: (i) they have the same number of elements; (ii) $\bar{B} \nsubseteq  B$; (iii)    $n/2$ elements of $\bar{B}$, namely, $n+2, n+ 4, \ldots , 2n-4, 2n-2, 2n,$ do not belong to $B$.

Suppose now that we wish to   multiply each element of the set of
natural numbers,~$\mathbb{N}$, by 2. Clearly, as the result we obtain the set, $\mathbb{E}$, of even numbers. Let us see whether the properties (i)--(iii) of the sets $B$ and $\bar{B}$ hold for $\mathbb{N}$ and~$\mathbb{E}$. With respect to  the property (i), we should say that, due to (\ref{even_bij}), the set obtained after multiplication has the same cardinality as the original set, i.e., it is countable. Then we see a paradoxical situation because,  in contrast with the   finite sets $B$ and $\bar{B}$, the set $\mathbb{E}$  obtained after multiplication is a \emph{proper subset} of  the original set $\mathbb{N}$, i.e., the property (ii) does not hold. Once again due to (\ref{even_bij}), the property (iii) does not hold either.

\subsection{Hilbert's paradox of the Grand Hotel}  This paradox proposed by David Hilbert in 1924 became popular thanks to the book ``One, Two, Three, ... Infinity'' of George Gamow (see \cite{Gamow}). It has the following formulation.  We all know that in a  hotel having a finite number   of rooms no more new
guests can be accommodated if it is full. Hilbert's Grand Hotel has
an infinite number of rooms (of course, the number of rooms is
countable, because the rooms in the Hotel are numbered). If a new
guest arrives at the Hotel where every room is occupied, it is,
nevertheless, possible to find a room for the newcomer. To do so, it
is necessary to move the guest occupying room 1 to room 2, the guest
occupying room 2 to room 3, etc. In such  way room 1 will be ready
for the new guest and, in spite of the assumption that there are no
available rooms in the Hotel, an empty room is found.

The paradox consists in the fact that we have supposed that the hotel is full and, nevertheless, it becomes possible to accommodate a newcomer in it. There exists different generalizations of this paradox showing in a similar way how it is possible to accommodate a finite and even infinite number of new guests in it.

\subsection{Three paradoxes regarding divergent series}

Let us now present another kind of paradoxes dealing with divergent series. We shall show that a very simple chain of equalities including addition of an infinite number of summands can lead
to a paradox. The first paradox is the following. Suppose that we have
\beq
x=1+2+4+8+\ldots
 \label{3.100}
       \eeq
Then we can multiple both parts of this equality by 2:
\[
2x=2+4+8+\ldots
\]
By adding 1 to both parts of the previous formula we obtain
 \beq
2x+1=1+2+4+8+\ldots
 \label{2.11}
       \eeq
It can be immediately noticed that the right hand side of
(\ref{2.11}) is just equal to $x$ and, therefore, it follows
\[
2x+1=x
\]
from which we obtain
\[
x=-1
\]
and, as a final paradoxical result, the following equality follows
\beq
 1+2+4+8+\ldots = -1.
\label{-1}
 \eeq
The paradox here is evident: we have summed up an infinite number of positive integers and have obtained as the final result a negative number.

The second paradox  considers    the well known divergent series of Guido Grandi $S = 1-1+1-1+1-1+\ldots$ By applying the telescoping rule, i.e., by writing a general element of the series as a difference, we can obtain two different answers    using two general elements, 1-1 and -1+1:
$$ S = (1-1)+(1-1)+(1-1)+\ldots =0,$$
$$ S = 1+(-1+1)+(-1+1) +(-1+1)+\ldots = 1.$$
In the literature there exist  many other approaches giving different answers regarding the value of this series (see, e.g., \cite{Knopp}). Some of them use various notions of
average (for instance, Ces\`aro summation assigns the value $0.5$ to $S$).

The third series we consider is the famous paradoxical result of Ramanujan
 \beq
c=1+2+3+ 4+ 5+ \ldots   = - 1/12.
 \label{Ramanu_3}
 \eeq
To obtain this remarkable result he multiplies
the left-hand part of (\ref{Ramanu_3}) by 4 and then subtracts the
result from (\ref{Ramanu_3}) as follows
 \beq
 \begin{array}{rclrrrrr}
  c & = & 1& + 2& + 3& + 4&+ 5 & + 6 +\ldots   \\
   4c & = & &  4& & + 8 & &+ 12+\ldots    \\
  -3c & = & 1& -2&+ 3&  -4&+  5  & -6+\ldots
   \end{array}
 \label{Ramanu_5}
 \eeq
Ramanujan then uses the result (considered in various forms by
Euler, Ces\`aro, and H\"older) attributing to the alternating series
$1 - 2 + 3 - 4 + \ldots$ the value $\frac{1}{4}$ as the formal power
series expansion of the function $\frac{1}{(1 + x)^2}$ for $x=1$,
i.e.,
 \beq
1 - 2 + 3 - 4 + \ldots = 1/4. \label{Ramanu_4}
 \eeq
 Thus, it follows from (\ref{Ramanu_5}) and (\ref{Ramanu_4}) that
\[
-3c  =1-2+ 3  -4+  5   -6+\ldots = 1/4,
 \]
from where Ramanujan gets (\ref{Ramanu_3}). This result looks even stranger than (\ref{-1}) because the sum of infinitely many positive integers is not only negative but also fractional.

\subsection{The rectangle paradox of Torricelli}  
\begin{figure}[t]
   \centering
   \includegraphics[width=12.5cm, height=8cm]{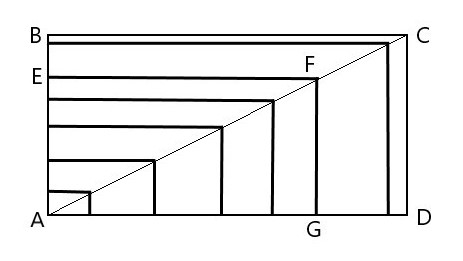}
  \caption{ The rectangle paradox of Torricelli  }
   \label{fig_1}
  \end{figure}

This paradox proposed by Evangelista Torricelli  (see, e.g. \cite{Alexander,Nillsen}) considers   a rectangle ABCD that is not a square (see Fig.~\ref{fig_1}). Without loss of generality let us suppose that the length $|AB|$ is two times smaller than $|BC|$. On the one hand, it is evident that the diagonal AC splits the rectangle   into two triangles ABC and CDA having equal areas. On the other hand, it is possible to propose the following reasoning using infinitesimals (Torricelli talks about indivisibles) that challenges this conclusion.

Let us  cover the
upper triangle ABC by an infinite number of horizontal line segments having an infinitesimal width. Analogously, the
lower triangle CDA is covered by the corresponding equal number  of vertical line segments also having an infinitesimal width. Fig.~\ref{fig_1}
illustrates only six   horizontal segments of this kind and the   corresponding six vertical segments.
By the construction, each horizontal line segment is two times greater in length than the corresponding vertical line
segment, for instance, $|EF|=2|FG|$. The area  of the upper triangle ABC can be obtained by summing up the areas of the horizontal line segments covering it. Analogously, the area  of the triangle ACD can be obtained by summing up the areas of the vertical lines covering it. Because each horizontal line segment has the length that is two times greater than the length of the corresponding vertical line segment, it   follows that the   triangle ABC has a greater area than the  triangle ACD. Thus, the paradox arises because   the two triangles have areas that are equal and in the same time are not equal.

\subsection{Thomson's lamp paradox}

The last paradox  we consider here is the
\textit{Thomson Lamp Paradox}
described in \cite{Thomson}. Suppose that we have a  lamp and
start to turn it on for $\frac{1}{2}$ minute, then turn it off for
$\frac{1}{4}$ minute, then on again for $\frac{1}{8}$ minute, etc.
At the end of one minute, the lamp switch will have been moved
infinitely many times (to be precise,   countably many times). Will then the lamp
be \textit{on} or \textit{off} at the end of one minute? It is easy to see that this
paradox is  equivalent to the following question: Is the `last'
integer  even or odd?

\section{Numeral systems of Pirah\~{a} and Munduruk\'u prompt a new point of view on infinity}
In order to understand how one can change his/her view on infinity,
let us consider some numeral systems used to express finite numbers.  Recall  that a \textit{numeral} is a symbol (or a
group of symbols) that represents a \textit{number} that is a
concept. The same number can be represented by different numerals.
For example, symbols `10', `ten', `IIIIIIIIII', and `X'  are different numerals, but they all represent the same
number. Rules used
to write down numerals together with algorithms for executing
arithmetical operations form a \textit{numeral system}. Thus,
numbers   can be considered as objects of an observation  that are
represented (observed)   by instruments of the observation, i.e., by
numerals and, more general, by numeral systems.

People in different historical periods used different numeral systems to count and these systems: (a) can be more or less suitable for counting; (b) can express different sets of numbers. For instance, Roman numeral system is not able to express zero
and negative numbers and such expressions as II -- VII or X -- X
are indeterminate forms in this numeral system. As a result, before
the appearance of positional numeral  systems and the invention of zero   mathematicians were not able to create
theorems involving zero and negative numbers and to execute
computations with them. The   positional numeral
system not only has allowed people to execute new operations but has
led two new theoretical results, as well. Thus, numeral systems not
only limit us in practical computations, they induce boundaries on
theoretical results, as well.

It should be stressed that the powerful positional numeral system also has its limitations. For example,  nobody is
able to write down a numeral in the decimal positional system having
$10^{100}$ digits  (see a discussion on feasible numbers in \cite{Parikh,Sazonov}). In fact, suppose that one is able to write down
one digit in one nanosecond. Then, it will take $10^{91}$ seconds to
record  all $10^{100}$ digits. Since in one year there are 31.556.926 $\approx 3.2 \cdot 10^{7}$ seconds, $    10^{91}$ seconds are
approximately $ 3.2 \cdot 10^{83}$ years. This is  a sufficiently
long time since it is supposed that the age of the universe is
approximately $  1.382 \cdot 10^{10}$ years.

As we have seen above,    Roman numeral system is weaker than the positional one. However,
it is not the weakest numeral system. There exist very poor numeral
systems allowing their users to express very few numbers and one of
them is illuminating for our story. This numeral system  is used by
a tribe, Pirah\~{a}, living in Amazonia nowadays. A study published
in \textit{Science} in 2004 (see \cite{Gordon}) describes that these
people use  an extremely simple numeral system for counting: one,
two, many. For Pirah\~{a}, all quantities larger than two are just
`many' and such operations as 2+2 and 2+1 give the same result,
i.e., `many'. Using their weak numeral system Pirah\~{a} are not
able to see, for instance, numbers 3, 4, and 5, to execute
arithmetical operations with them, and, in general, to say anything
about these numbers because in their language there are neither
words nor concepts for that.

It is worthy of mention that the result `many' is not
wrong\label{pr:4}. It is just \textit{inaccurate}.  Analogously,
when we observe a garden with 547 trees, then both phrases: `There
are 547 trees in the garden' and `There are many trees in the
garden' are  correct. However, the accuracy of the former phrase is
higher than the accuracy of the latter one. Thus, the introduction
of a numeral system having numerals for expressing numbers 3 and 4
leads to a higher accuracy of computations and allows one to
distinguish results of operations 2+1 and 2+2.

The poverty of the  numeral system of Pirah\~{a} leads   also to the
following results
 \beq \mbox{`many'}+ 1= \mbox{`many'},
\hspace{2mm}    \mbox{`many'} + 2 = \mbox{`many'},\hspace{2mm}
\mbox{`many'}+ \mbox{`many'} = \mbox{`many'} \label{1}
 \eeq
that are crucial for changing   our outlook on infinity. In fact, by
changing in these relations `many' with $\infty$ we get relations
used to work with infinity   in the traditional calculus
 \beq
\infty + 1= \infty,    \hspace{5mm}    \infty + 2 = \infty,
\hspace{5mm}   \infty +
  \infty =   \infty.
  \label{2}
 \eeq

Analogously, if we consider Cantor's cardinals (where, as usual,
numeral $\aleph_0$ is used for cardinality of countable sets and
numeral
 \textfrak{c}  for cardinality of the continuum, see, e.g.,
\cite{Woodin}) we
have similar relations
 \beq
 \aleph_0+  1 = \aleph_0, \hspace{5mm}    \aleph_0 + 2 =\aleph_0,
 \hspace{5mm}  \aleph_0+ \aleph_0 = \aleph_0,
  \label{3}
 \eeq
 \beq
 \textfrak{c} +  1 = \textfrak{c} , \hspace{5mm}    \textfrak{c}  + 2 = \textfrak{c} ,
 \hspace{5mm}  \textfrak{c}  + \textfrak{c}  = \textfrak{c} .
 \label{4}
 \eeq

It should be mentioned that the astonishing numeral system of
Pirah\~{a} is not an isolated example of this way of counting. In \cite{Comrie}, more than 20 languages having numerals only for small numbers are mentioned. For example, the same counting system, one, two, many, is used by the
Warlpiri people, aborigines living in  the Northern Territory of
Australia (see \cite{Butterworth}).   The Pitjantjatjara people
living in the Central Australian desert use  numerals    one, two,
three, big mob (see \cite{Leder}) where `big mob' works as `many'.
It makes sense to remind also another Amazonian tribe -- Munduruk\'u
(see \cite{Pica}) who fail in exact arithmetic with numbers larger
than 5 but are able to compare and add large approximate numbers
that are far beyond their naming range. Particularly, they use the
words `some, not many' and `many, really many' to distinguish two
types of large numbers. Their arithmetic with  `some, not many' and
`many, really many' reminds    the rules Cantor uses to
work with $\aleph_0$ and \textfrak{c}, respectively. In fact, it is
sufficient to compare
 \beq
\mbox{`some, not many'+ `many, really many' =   `many, really many'}
  \label{5}
 \eeq
 with
 \beq
 \aleph_0+ \textfrak{c}  = \textfrak{c}  \label{6}
 \eeq
to see this similarity.

Let us compare now the weak numeral systems involved in (\ref{1}),
(\ref{5}) and numeral systems used to work with infinity. We have
already seen that relations (\ref{1}) are results of the weakness of
the   numeral system employed. Moreover, the usage of a stronger
numeral system shows that  it is possible to pass  from records
1+2~=~`many' and 2+2~=~`many' providing for two different
expressions the same result, i.e., `many',  to more precise answers
1+2~=~3 and 2+2~=~4 and to see that $3 \neq 4$. In these examples we
have the same objects  -- small finite numbers -- but results of
computations we execute are different in dependence of the
instrument -- numeral system  -- used to represent numbers.
Substitution of the numeral `many' by a variety of numerals
representing numbers 3, 4, etc. allows us both to avoid relations of
the type  (\ref{1}), (\ref{5}) and to increase the accuracy of
computations.

Relations (\ref{2})--(\ref{4}), (\ref{6})   manifest a complete
analogy with (\ref{1}), (\ref{5}).   Canonically, symbols $\infty,$
$\aleph_0,$ and \textfrak{c}  are \textit{identified} with concrete
mathematical objects and (\ref{2})--(\ref{4}), (\ref{6}) are
considered as intrinsic properties of these infinite objects (see
e.g.,
\cite{Robinson,Woodin}).
However, the analogy with (\ref{1}), (\ref{5}) suggests that relations
(\ref{2})--(\ref{4}), (\ref{6}) do not reflect the \textit{nature}
of   infinite objects. They are just  a result of weak numeral
systems used to express infinite quantities. As (\ref{1}), (\ref{5})
show the lack of numerals in numeral systems of Pirah\~{a},
 Warlpiri, Pitjantjatjara, and  Munduruk\'u for expressing different
finite quantities, relations (\ref{2})--(\ref{4}), (\ref{6}) show
shortage of numerals in mathematical analysis   and in set
theory for expressing different infinite numbers. Another hint
leading to the same conclusion is the situation with indeterminate
forms of the kind III-V in Roman numerals that have been excluded
from the practice of computations after introducing positional
numeral systems.

Thus, the analysis made above allows us to formulate the following
key observation that changes our perception of infinity:
\begin{quotation}
  \textbf{Our difficulty in working with infinity is not a
consequence of the \textit{nature} of infinity but is a result of
\emph{weak numeral systems} having too little numerals  to express
the multitude of infinite numbers.}
 \end{quotation}

The way of reasoning where the object of the study is separated from
the tool used by the investigator  is very common in natural
sciences where researchers use tools to describe the object of their
study and the used instrument influences the results of the
observations and determine their accuracy. The same happens in
Mathematics studying natural phenomena, numbers,
 objects that can be constructed by using numbers, sets, etc. Numeral
systems used to express numbers are among the instruments of
observation  used by mathematicians. As we have illustrated above,
the usage of powerful numeral systems gives the possibility to
obtain more precise results in Mathematics in the same way as usage
of a good microscope gives the possibility of obtaining more precise
results in Physics. Traditional numeral systems have been developed
to express finite quantities and they simply have no sufficiently
high number of numerals to express different infinities (and
infinitesimals).

\section{A new way of counting}
 In this section, we briefly describe a recent  numeral system that can be used to write down various infinite, finite, and infinitesimal numbers in a unique framework (see a comprehensive technical survey \cite{EMS}  and a popular book \cite{Sergeyev} for its description) concentrating ourselves on details that then will be used to reconsider the paradoxes. It should be emphasized immediately that the  methodology to be
presented     is not a
contraposition to the ideas of Cantor, Levi-Civita, and Robinson. In
contrast, it is   an applied evolution of their ideas.
The new computational methodology introduces the notion of the
accuracy of numeral systems and shows that different numeral systems
can express different sets of finite and infinite numbers with
different accuracies. A clear analogy with Physics can be established in this context.

When a physicist uses a weak lens $A$
and sees two black dots in his/her microscope he/she does not say:
The object of the observation \textit{is} two black dots. The
physicist is obliged to say: the lens used in the microscope allows
us to see two black dots and it is not possible to say anything more
about the nature of the object of the observation until we replace
the instrument - the lens or the microscope itself - with a more
precise one. Suppose that  he/she changes the lens and uses a
stronger lens $B$ and is able to observe that the object of the
observation is viewed as ten (smaller) black dots. Thus, we have two
different answers: (i) the object is viewed as two dots if the lens
$A$ is used; (ii) the object is viewed as ten dots by applying the
lens $B$. Which of the answers is correct? Both. Both answers are
correct but with the \textit{different accuracies} that depend on
the lens used for the observation. The answers are not in opposition
one to another, they both describe the reality (or whatever is
behind the microscope) correctly with the precision of the used
lens. In both cases our physicist discusses what he/she observes and
does not pretend to say what the object \emph{is}.

We shall do the same with infinite numbers and sets (objects of our study) and numeral systems used to observe them (our tools).  Traditional
approaches and the methodology described here do not contradict one
another, they are just different lenses having different
accuracies for observations of mathematical objects.

Before we start a technical consideration let us mention that a number of papers studying consistency of the new methodology  and its connections   to the historical panorama of ideas dealing with infinities and
infinitesimals have been
published (see
\cite{Gangle,Lolli,Lolli_2,MM_bijection,Sorbi,first,Fallacies,Tohme2020}). In particular, in \cite{Fallacies} it is stressed that it is not related to non-standard analysis.   The
methodology   has   been   successfully applied in
several areas of mathematics and computer science (more than 60 papers   published in international scientific journals can be found at the dedicated web page \cite{www}). We provide here just a few examples of areas where this methodology is useful. First of all, its successful applications in teaching mathematics should be mentioned   (see, e.g., \cite{Antoniotti_teaching,Iannone,Ingarozza_teaching} and the dedicated web page \cite{web_teaching} developed at the University of East Anglia, UK and containing, among other things, a comprehensive teaching manual and a nice animation related to the Hilbert's paradox of   the Grand Hotel). Then, we can indicate game theory and probability (see, e.g., \cite{Calude_Dumitrescu,DAlotto_games,Fiaschi_Cococcioni_2018,Fiaschi2020,Pepelyshev_Zhigljavsky,Rizza_2,Rizza_3}); local, global, and
multiple criteria optimization  (see
\cite{Cococcioni,DeLeone,DeLeone_2,Fasano,Gaudioso&Giallombardo&Mukhametzhanov(2018),medals,homogeneity,Zilinskas}), hyperbolic
geometry   and percolation (see
\cite{Iudin,Margenstern}), fractals (see
\cite{Antoniotti2020,Caldarola_1,Koch}),  infinite
series   (see
\cite{EMS,Zhigljavsky}),
  Turing machines, cellular automata,  and supertasks  (see
\cite{DAlotto,Rizza,Rizza_3,Sergeyev_Garro}),   numerical
differentiation and numerical solution of  ordinary differential
equations (see \cite{ODE_3,Falcone:et:al.2020,Falcone:et:al_Zeno_1,ODE_5}), etc.

In order to start, let us mention that for thousands of years on the Earth there exists a way of counting   huge finite quantities that has not been formalized until the recent times. In traditional mathematics, after appearance  of axioms of Peano, natural numbers are introduced starting from 0 by adding a unit to  get 1, then adding another unit to 1 to obtain 2, and, by continuing in this way, other positive integers are introduced. Let us illustrate   by an example  that counting  is a more complex procedure with respect to just adding 1 to 0 many times.

Imagine that we are in a granary and the owner asks us to count how
much grain he has inside it. Obviously, it is possible to answer
that there are many seeds in the granary. This answer is correct but
its accuracy is low. In order to obtain a more precise answer it
would be necessary to count the grain seed by seed but since the
granary is huge, it is not possible to do this due to practical
reasons.

To overcome this difficulty and to obtain an answer that is more accurate than
`many', people take sacks, fill them   with seeds, and count the
number of sacks. In this situation, we suppose that: (i) all the
seeds have the same measure and all the sacks also; (ii) the number
of seeds in each sack is the same and is equal to $K_1$ but the sack
is so big  that we are not able to count how many seeds it contains
and to establish the value of $K_1$; (iii) in any case the resulting
number $K_1$ would not be expressible by available numerals.

Then, if the granary is huge and it becomes difficult to count the
sacks, then trucks or even big train wagons are used. As it was for
the sacks, we suppose   that all trucks contain the same number
$K_2$ of sacks, and all train wagons contain the same number $K_3$
of trucks, however, the  numbers $K_i, i=1,2,3,$ are so huge that it
becomes impossible to determine their values.
  At the end of   this counting   we
obtain a result in the following form: the granary contains 34
wagons, 27 trucks, 16 sacks, and 134 seeds of grain. Note, that if
we add, for example, one seed to the granary, we can count it and
not only see that  the granary has more grain but also quantify the
increment: from 134 seeds we pass to 135 seeds. If we take out one wagon, we again are able to say how
much grain has been subtracted: from 34 wagons we pass to 33 wagons.

Let us make some considerations upon the way of counting described
above. In our example it is necessary to count large quantities.
They are finite but it is impossible to count them directly by using
the elementary unit  of measure, $u_0$, (seeds) because the
quantities expressed in these units would be too large. Therefore,
people are forced to behave as if the quantities were infinite.

To solve the problem of `infinite' quantities, new units of measure,
 $u_1$ -- sacks, $u_2$ -- trucks, and $u_3$ -- wagons, are
introduced. The new units have   an important feature: all the units
$u_{i+1}$ contain a certain number $K_i$\label{pr:2} of units
$u_{i}$ but these numbers, $K_i, i=1,2,3,$ \emph{are unknown}.
  Thus, quantities that it was impossible to
express using only the initial unit  of measure, $u_0$, are
perfectly expressible in the new units  $u_i, i=1,2,3$. Notice that, in spite of the fact that the numbers $K_i$ are
unknown, the accuracy of the obtained answer is equal to one seed.
In fact, if we add one seed we are able to register and to quantify
that we have more seeds and if we subtract one wagon and  two
sacks, we again can quantify the decrease.

 This key idea of counting by introduction of
new units of measure with unknown but fixed values $K_i$   will be
used in what follows to deal with infinite quantities together with
the    relaxation allowing one to use negative digits in positional
numeral systems. It is necessary to extend the idea of the
introduction of new units of measure from sets and numbers that are
huge but finite to infinite sets and numbers. This can be done by
extrapolating from finite to infinite the idea that $n$ is both the
number of elements of the set $\{ 1, 2, 3, \ldots , n-1, n \}$ and the last element of this set.

The infinite unit of measure is  introduced as the number of
elements of  the set,~$\mathbb{N}$, of natural
numbers  and expressed by the numeral
\ding{172} called \textit{grossone}. Using the granary example
discussed above we can offer the following interpretation: the set
$\mathbb{N}$ can be considered as a sack and \G1 is the number of
seeds in the sack. Following our extrapolation, the introduction of
\G1 allows us to write down the set of natural numbers in the form
 \beq
 \mathbb{N} = \{ 1, 2, 3, \ldots , \G1-3, \G1-2, \G1-1, \G1 \},
  \label{N}
       \eeq
  where  $\G1-3, \G1-2, \G1-1, \G1$ are infinite natural numbers. Thus, the set of natural numbers will be written in the form (\ref{N}) instead of the usual record (\ref{paradoxes_N}). We emphasize that in both cases we deal
with the same mathematical object -- the set of natural numbers --
that is observed through two different instruments. In the traditional
case, usual numeral systems do not allow us to express infinite
numbers   whereas the numeral system with grossone offers this
possibility. Similarly,   Pirah\~{a}  are not able to see finite
natural numbers greater than~2 but these numbers (e.g., 3 and 4)
belong to $\mathbb{N}$ and are visible if one uses a more powerful
numeral system. Notice also that in traditional statements (for example, in non-standard analysis) infinite numbers are not included in $\mathbb{N}$. However, if it is supposed that $\mathbb{N}$ is infinite and its elements are constructed, starting from 1 (or zero, as Peano did), according to the rule: the number $n$ is followed by the number $n + 1$, then each next number will be finite and, therefore, all natural numbers will be finite. Thus, by this construction, any set $\{1, 2, 3, ..., n \}$ will contain a finite number of elements. This would contradict  the assumption that $\mathbb{N}$ is an infinite set.

Grossone  is introduced by describing its properties  postulated by
the \textit{Infinite Unit Axiom} (IUA) consisting of three parts:
Infinity, Identity, and Divisibility. Similarly, in order to pass
from natural to integer numbers a new element -- zero -- is
introduced, a numeral to express it is chosen, and   its properties
are described. The IUA is added to axioms for real numbers. Thus, it is postulated that
associative and commutative properties of multiplication and
addition, distributive property of multiplication over addition,
existence of   inverse elements with respect to addition and
multiplication hold for grossone as they do for finite
numbers.

Let us introduce the axiom and then give some comments upon it.
Notice that in the IUA infinite sets will be described in the
traditional form, i.e., without indicating the last element. For
instance, the set of natural numbers will be written as   (\ref{paradoxes_N})  instead of the record  (\ref{N}) that will be used after the axiom will be introduced.

\textbf{The  Infinite Unit Axiom.} The infinite unit of measure  is
introduced as the number of elements of  the set, $\mathbb{N}$, of
natural numbers. It is expressed by the numeral \ding{172} called
\textit{grossone} and has the following properties:

\textit{Infinity.}
  Any finite  natural number $n$  is less than grossone, i.e.,  $n
<~\G1$.

\textit{Identity.}
 The following
relations  link \ding{172} to identity elements 0 and 1
 \beq
 0 \cdot \G1 =
\G1 \cdot 0 = 0, \hspace{3mm} \G1-\G1= 0,\hspace{3mm}
\frac{\G1}{\G1}=1, \hspace{3mm} \G1^0=1, \hspace{3mm}
1^{\mbox{\tiny{\ding{172}}}}=1, \hspace{3mm}
0^{\mbox{\tiny{\ding{172}}}}=0.
 \label{3.2.1}
       \eeq

\textit{Divisibility.} For any finite natural number  $n$   sets
$\mathbb{N}_{k,n}, 1 \le k \le n,$ being the $n$th parts of the set,
$\mathbb{N}$, of natural numbers have the same number of elements
indicated by the numeral $\frac{\G1}{n}$ where
 \beq
 \mathbb{N}_{k,n} = \{k,
k+n, k+2n, k+3n, \ldots \}, \hspace{5mm} 1 \le k \le n,
\hspace{5mm} \bigcup_{k=1}^{n}\mathbb{N}_{k,n}=\mathbb{N}.
 \label{3.3}
       \eeq

Let us comment upon this axiom. Its first part   -- Infinity -- is
quite clear. In fact, we want to describe an infinite number, thus,
it should be larger than any finite number. The second part of the
axiom -- Identity -- tells us that \ding{172} interacts  with
identity elements~0 and 1 as all other numbers do. In   the
moment when we have stated that grossone is a number, we have fixed the
usual properties of numbers, i.e., the properties described in
Identity, associative and commutative properties of multiplication
and addition, distributive property of multiplication over addition,
etc. The third part of the axiom -- Divisibility -- is the most
interesting, since it links infinite numbers to infinite sets (in
many traditional theories  infinite numbers are introduced
algebraically, without any connection to infinite sets). It  is
based on Euclid's   Common Notion no.~5 `The whole is greater than the part'. In the new methodology, it is applied to all quantities: finite, infinite, and infinitesimals.

Let us consider   two examples for $n=1$ and $n=2$ in (\ref{3.3}).   If we take $n = 1$, then it follows that
$\mathbb{N}_{1,1} = \mathbb{N}$ and Divisibility says that the set,
$\mathbb{N}$, of natural numbers has \ding{172} elements. If $n =
2$, we have two sets $\mathbb{N}_{1,2}$ and $\mathbb{N}_{2,2}$,
where
 \beq
\begin{array}{ccccccccccc}
     \mathbb{N}_{1,2} = &  \{1, &   & 3, &  & 5, &   & 7, & \ldots  &\}, \\
      &   &   &   &   &   &   &    &  \\
     \mathbb{N}_{2,2} = &  \hspace{-3mm}  \{  & 2, &  & 4, &  & 6, &   &  \ldots &\}  \\
    \end{array}
\label{3.3.1}
       \eeq
and they have $\frac{\G1}{2}$  elements each. Notice that the sets
$\mathbb{N}_{1,2}$ and $\mathbb{N}_{2,2}$ have the same number of
elements not because they are in a one-to-one correspondence but due to
the Divisibility axiom.  In fact, we are not able to count the number
of elements of the sets $\mathbb{N}$, $\mathbb{N}_{1,2}$, and
$\mathbb{N}_{2,2}$ one by one because   we are able
to execute only a finite number of operations (we emphasize here the practical orientation of this methodology) whereas   these sets are
infinite. To define their number of elements we use Divisibility and
implement the principle `The whole is greater than the part' in practice by determine the number of the elements of the parts using the whole.

In general, to introduce $\frac{\G1}{n}$ we do not try to count
elements $k, k+n, k+2n, k+3n, \ldots$ one by one in (\ref{3.3}). In
fact, we cannot do this due to the finiteness of our practical counting abilities. By using Euclid's principle, we
construct the sets $\mathbb{N}_{k,n}, 1 \le k \le n,$ by separating
the whole, i.e., the set $\mathbb{N}$, in $n$ parts  and   we affirm
that the number of elements of the $n$th part of the set, i.e.,
$\frac{\G1}{n}$, is $n$ times less than the number of elements of
the entire set, i.e., than \ding{172}.

As was already mentioned, in terms of our granary example \ding{172}
can be interpreted as the number of seeds in the sack. In that
example, the number $K_1$ of seeds in each sack was fixed and finite
but it was impossible to  express it in units $u_0$, i.e., seeds, by
counting seed by seed because we had supposed that sacks were very
big and the corresponding number would not be expressible by
available numerals. In spite of the fact that $K_1, K_2,$ and $ K_3$
were inexpressible and unknown,  by using new units of measure
(sacks, trucks, etc.) it was possible to count more easily  and to
express the required quantities. Now our sack has the infinite but
again \textit{fixed} number of seeds. It is fixed because it has a
strong link to a concrete set -- it is the number of elements of the
set of natural numbers. Since this number is inexpressible by
existing numeral systems with the same  accuracy afforded to measure
finite small sets, we introduce a
new numeral, \G1, to express  the required quantity. Then, we apply
Euclid's principle and say that if the sack contains~\ding{172} seeds,
then, even though we are not able to count the number of seeds of
the $n$th part of the sack seed by seed, its $n$th part contains $n$
times less seeds than the entire sack, i.e., $\frac{\G1}{n}$ seeds.
Notice that    the numbers $\frac{\G1}{n}$ are integer since they
have been introduced as numbers of elements of sets
$\mathbb{N}_{k,n}$.

The new unit of measure allows us to express a variety of infinite
numbers (including those larger than \G1  that will be considered
shortly) and calculate easily the number of elements of the union,
intersection, difference, or product of   sets of   type
$\mathbb{N}_{k,n}$. Due to our accepted methodology, we do it in the
same way  as these measurements are executed for finite sets. Let us
consider two simple examples showing how grossone can be used for
this purpose (see \cite{EMS} for a detailed discussion).

    Let us determine the number of elements of the set
$A_{k,n} = \mathbb{N}_{k,n} \backslash \{a\},$ $a \in
\mathbb{N}_{k,n}, n \ge 1$. Due to the IUA, the set
$\mathbb{N}_{k,n}$ has $\frac{\G1}{n}$ elements. The set $A_{k,n}$
has been constructed by excluding one element from $N_{k,n}$. Thus,
the set $A_{k,n}$ has $\frac{\G1}{n}-1$ elements. The granary
interpretation can be also given for the number $\frac{\G1}{n}-1$ as
the number of seeds in the $n$th part of the sack minus one seed.
For $n=1$ we have $\G1-1$ interpreted as the number of seeds in the
sack minus one seed.

 Let us consider the following two sets
\[
B_1 = \{ 4, 9, 14, 19, 24, 29, 34, 39, 44, 49, 54, 59, 64, 69, 74,
79,\ldots\},
\]
\[
B_2 = \{ 3, 14, 25, 36, 47, 58, 69, 80, 91, 102, 113, 124, 135,
\ldots\}
\]
and determine the number of elements   in the set $B  = (B_1 \cap
B_2 ) \cup \{  3,4,5, 69  \}$. It follows immediately from the IUA
that $B_1 = \mathbb{N}_{4,5}$ and $B_2 = \mathbb{N}_{3,11}$. Their
intersection
\[
B_1 \cap B_2 = \mathbb{N}_{4,5} \cap \mathbb{N}_{3,11} = \{ 14,
69, 124, \ldots\} = \mathbb{N}_{14,55}
\]
and, therefore, due to the IUA, it has $\frac{\G1}{55}$ elements.
Finally, since 69 belongs to the set  $\mathbb{N}_{14,55}$ and 3, 4,
and 5 do not belong to it, the set $B$ has $\frac{\G1}{55}+3$
elements. The granary interpretation: this is the number of seeds in
the $55$th part of the sack plus three seeds.

The IUA introduces \G1 as the number of elements of the set of natural
numbers and, therefore, it is the last natural number. We can also talk about the set of \textit{extended natural
numbers} indicated as $\widehat{\mathbb{N}}$ and including
$\mathbb{N}$ as a proper subset
\[
  \widehat{\mathbb{N}} = \{
\underbrace{1,2, \ldots ,\mbox{\ding{172}}-1,
\mbox{\ding{172}}}_{\mbox{{\scriptsize Natural numbers}}},
\mbox{\ding{172}}+1, \mbox{\ding{172}}+2, \ldots ,
2\mbox{\ding{172}} -1, 2\mbox{\ding{172}}, 2\mbox{\ding{172}} +1,
\ldots
\]
 \beq
\mbox{\ding{172}}^2-1, \mbox{\ding{172}}^2, \mbox{\ding{172}}^2+1,
\ldots 3\mbox{\ding{172}}^{\mbox{\tiny{\ding{172}}}}-1,
3\mbox{\ding{172}}^{\mbox{\tiny{\ding{172}}}},
3\mbox{\ding{172}}^{\mbox{\tiny{\ding{172}}}}+1, \ldots \}.
\label{4.2.2}
       \eeq
The extended natural numbers  greater than grossone are also
linked to infinite sets of numbers and can   be interpreted in the terms of
grain. For example, \G1+1 is the number of elements of a set $B_3 = \mathbb{N}   \cup \{a\},$ where $a$ is integer and $a \notin \mathbb{N}$. In the terms of grain, \G1+1 is the number of seeds in a sack plus one seed.

   Let us give another example and determine the number of elements of the set
 \[
 B_4  =
\{
  (a_1, a_2) : a_i \in   \mathbb{N},   i \in   \{1,2\}
 \},
 \]
being the set of couples of natural numbers.  It is
known from combinatorial calculus that if we have two positions and
each of them can be filled in by one of $l$ symbols, the number of
the obtained couples is equal to $l^2$. In our case, since
$\mathbb{N}$  has grossone elements, $l = \G1$. Thus, the set $B_4$
has  $\G1^2$ elements.
This fact is illustrated below
\[
  \begin{array}{ccccc}
(1,1), & (1,2), & \ldots & (1,\G1-1), & (1,\G1), \\
  (2,1), & (2,2), & \ldots & (2,\G1-1), & (2,\G1), \\
\ldots & \ldots & \ldots & \ldots & \ldots \\
\hspace{3mm}(\G1-1,1), \hspace{3mm}& \hspace{3mm} (\G1-1,2),\hspace{3mm} & \hspace{3mm} \ldots \hspace{3mm} & \hspace{3mm}(\G1-1,\G1-1), \hspace{3mm} & (\G1-1,\G1), \\
(\G1,1), & (\G1,2), & \ldots & (\G1,\G1-1), & (\G1,\G1).
  \end{array}
\]

The introduced numeral system allows us to observe not only initial elements of certain infinite sets but also the final ones and some other infinite numbers in these sets. For example, we can write now the following records
\[
\mathbb{N} = \{1,2, \ldots , \ldots \frac{\mbox{\ding{172}}}{2}-1,  \frac{\mbox{\ding{172}}}{2}, \frac{\mbox{\ding{172}}}{2}+1, \ldots  \mbox{\ding{172}}-1,
\mbox{\ding{172}} \},
\]
\[
\mathbb{O} = \{1, 3, 5,  \ldots , \ldots \frac{\mbox{\ding{172}}}{2}-1,    \frac{\mbox{\ding{172}}}{2}+1, \ldots  \mbox{\ding{172}}-3,
\mbox{\ding{172}}-1 \},
\]
\[
\mathbb{E} = \{ 2, 4, 6, \ldots , \ldots \frac{\mbox{\ding{172}}}{2}-2,  \frac{\mbox{\ding{172}}}{2}, \frac{\mbox{\ding{172}}}{2}+2, \ldots  \mbox{\ding{172}}-2,
\mbox{\ding{172}} \},
\]
\[
\mathbb{Z} = \{-\G1, -\G1+1, -\G1+2, \ldots -2, -1, 0, 1,2, \ldots ,  \mbox{\ding{172}}-1,
\mbox{\ding{172}} \}.
\]
Due to the IUA, the set, $\mathbb{O}$, of odd numbers has $\frac{\mbox{\ding{172}}}{2}$ elements, the set, $\mathbb{E}$, of even numbers also has $\frac{\mbox{\ding{172}}}{2}$ elements. It is easy to calculate the number of elements of the set, $\mathbb{Z}$, of integers. It has \G1 positive elements, \G1 negative ones, and zero. Thus, the set  $\mathbb{Z}$ has $2\G1+1$ elements. For the purpose of this article the introduced material is sufficient. As was already mentioned, more information about \G1 can be found in a comprehensive technical survey \cite{EMS} and a popular book \cite{Sergeyev}.

In order to conclude this section it is necessary to emphasize that
the introduced numeral system  cannot give answers to \textit{all} questions regarding
infinite sets. As all numeral systems, it has its limitations. What can we say, for instance, about the number of
elements of the set  $\widehat{\mathbb{N}}$? Was this set described completely? The
introduced numeral system based on \ding{172} is too weak to give
answers to these questions since it does not allow us to express the
number of elements of this set. It is necessary to introduce in a reasonable
way a more powerful numeral system by defining new numerals (for
instance, \ding{173}, \ding{174}, etc).

\section{Paradoxes of infinity revisited} In this section we reconsider in the \G1-based  framework the paradoxes described above. It should be mentioned that several other  paradoxes related to infinities and infinitesimals in  probability theory and decision making we considered using the \G1-based methodology   in \cite{Fiaschi2020,Rizza,Rizza_2,Rizza_3}.

\subsection{Galileo's paradox}   Before we start to consider the set $I^2$  of square natural numbers from (\ref{4.3.0.1}), let us study bijection (\ref{even_bij}) between the sets of even and natural numbers. The traditional conclusion from (\ref{even_bij}) is that both sets
are countable and they have the same   cardinality~$\aleph_0$.

Let us see now what we can say from the new methodological position,
in particular, by using Euclid's principle together with the separation of the objects of study form the tools used for this study. The objects
of study here are  two infinite sets, $\mathbb{N}$ and $\mathbb{E}$, and the   instrument used to
compare them is the bijection. Since we know that some elements of $\mathbb{N}$ do not belong to $\mathbb{E}$, the separation of the objects of our study from the tools suggests that another
conclusion can be derived from (\ref{even_bij}): the accuracy of the
used instrument is not sufficiently high to see the difference
between the sizes of the two sets.

We have already seen that when one executes the operation of
counting, the accuracy of the result depends on the numeral system
used for counting. If one asked Pirah\~{a} to measure sets
consisting of four apples and five apples the answer would be that
both sets of apples have many elements. This answer is correct but
its precision is low due to the weakness of the numeral system used
to measure the sets.

 Thus, the introduction of the notion of accuracy for measuring sets  is
very important and should be applied to infinite sets also. As was
already discussed earlier, the similarity of Pirah\~{a}'s rules (\ref{1}) with the relations (\ref{3}) and (\ref{4}) holding for cardinal numbers  suggests
that the accuracy of the cardinal numeral
system of Alephs   is not sufficiently high  to see the difference
with respect to the number of elements of the two sets from
(\ref{even_bij}).

In order to look at the record (\ref{even_bij}) using the new
methodology, let us remind  that due to the IUA  the sets of even and
odd numbers have $\G1/2$ elements each and, therefore, \ding{172} is
even. It is also necessary to recall that numbers that are larger
than \G1 are not natural, they are  extended natural numbers. For
instance, $\G1+2 $ is even but not natural, it is the extended
natural, see (\ref{4.2.2}). Thus, the last even natural number is
\G1. Since the number of elements of the set of even numbers is
equal to $\frac{\G1}{2}$, we   can write down not only the initial
(as it is usually done traditionally) but also the final part of
(\ref{even_bij})
  \beq
\begin{array}{cccccccccc}
 2, & 4, & 6, & 8,  & 10, & 12, & \ldots  &
\G1 -4,  &    \G1  -2,   &    \G1    \\
 \updownarrow &  \updownarrow & \updownarrow  &
\updownarrow  & \updownarrow  &  \updownarrow  & &
  \updownarrow    & \updownarrow   &
  \updownarrow
   \\
 1, &  2, & 3, & 4 & 5, & 6,   &   \ldots  &    \frac{\G1}{2} - 2,  &
     \frac{\G1}{2} - 1,  &    \frac{\G1}{2}   \\
     \end{array}
\label{even_bij2}
 \eeq
concluding so (\ref{even_bij})   in a complete accordance with the
principle `The part is less than the whole'. Both records,
(\ref{even_bij}) and (\ref{even_bij2}), are correct but
(\ref{even_bij2}) is more accurate, since it allows us to observe the
final part of the correspondence that is invisible if
(\ref{even_bij}) is used. The new \G1-based numerals allow us to see that the set $\mathbb{E}$ has two times less elements with respect to  $\mathbb{N}$. Thus, the paradox is not present.

We are ready now to consider the set $I^2$  of square natural numbers from (\ref{4.3.0.1}) and the corresponding bijection from (\ref{trad_square}). The traditional reasoning does not allow one to see that these two sets have   different
numbers of elements. Again, the answer that both sets are countable is
correct but its accuracy is low. The \G1-based methodology allows us
to see the difference in their number of elements and to express the
final part of (\ref{trad_square}). The set $I^2$ can be defined now more accurately by emphasizing the fact that, by definition of \G1, each square natural number should be less than  or equal to grossone:
 \[
 I^2 = \{ x : x \in \mathbb{N}, i \in \mathbb{N}, x =i^{\,2}, x \le \G1 \}.
  \}
 \]
Then,  the number of elements, $J$, of the set $I^2$  can be determined  as
 \[
 J = \max
\{i : i^{\,2} \le \G1 \}. \label{I2}
 \]
By solving the required inequality $i^{\,2} \le \G1$ and taking the maximal integer $i$ satisfying we obtain that $J= \lfloor  \G1^{1/2}\rfloor$.   Thus, as it was with the sets $\mathbb{N}$ and $\mathbb{O}$, we can re-write the bijection (\ref{trad_square}) indicating both its initial and finite elements:
  \beq
 \begin{array}{ccccccccc}
 1, & 2^2, & 3^2, & 4^2,  & 5^2, &  \ldots  &
 (\lfloor  \G1^{1/2}\rfloor-2)^2 ,  &    (\lfloor  \G1^{1/2}\rfloor-1)^2,    &   \lfloor  \G1^{1/2}\rfloor^2   \\
 \updownarrow &  \updownarrow &
\updownarrow  & \updownarrow  &  \updownarrow  & &
  \updownarrow    & \updownarrow   &
  \updownarrow
   \\
 1, &  2, & 3, & 4 & 5, &    \ldots  &    \lfloor  \G1^{1/2}\rfloor   - 2,  &
      \lfloor  \G1^{1/2}\rfloor   - 1,  &     \lfloor  \G1^{1/2}\rfloor    \\
     \end{array}
\label{new_square}
 \eeq
 Since $\lfloor  \G1^{1/2}\rfloor < \G1$, this paradox also vanishes.

 \subsection{Set-multiplication paradox}

The introduction of the \G1-based numeral system allows us   to
write down the sets of natural  and extended natural numbers in
the form (\ref{4.2.2}).
 By definition, the number of elements of~$\mathbb{N}$ is equal to \G1. After multiplication of each of
the elements of $\mathbb{N}$  by 2,  the resulting set, that we call
$\mathbb{E}^2$,  will also have grossone elements because muliplication of elements of a set by a constant that is not equal to zero does not change the number of elements of the set. In fact,
the number $\frac{\G1}{2}$ multiplied by 2 gives us \G1 and
$\frac{\G1}{2}+1$ multiplied by 2 gives us $\G1+2 $ that is even
extended natural number, see (\ref{4.2.2}). Analogously, the last
element of $\mathbb{N}$, i.e., \G1, multiplied by 2 gives us 2\G1.
Thus, the set  $\mathbb{E}^2$ can be written as follows
\[
 \mathbb{E}^2 = \{  2, 4, 6, \hspace{2mm}  \ldots  \hspace{2mm}
\G1-2, \G1, \G1+2,\hspace{2mm}  \ldots \hspace{2mm} 2\G1-4, 2\G1-2,
  2\G1 \}
  \]
  and the corresponding bijection is
  \[
\begin{array}{ccccccccccc}
 2, & 4, & 6, & \ldots & \G1-2,& \G1, &\G1+2, & \ldots & 2\G1-4,  & 2\G1-2, & 2\G1   \\
 \updownarrow &  \updownarrow & \updownarrow  & &
\updownarrow  & \updownarrow  &  \updownarrow  & &
  \updownarrow    & \updownarrow   &
  \updownarrow
   \\
 1, &  2, & 3, & \ldots & \frac{\G1}{2}-1, & \frac{\G1}{2},   & \frac{\G1}{2}+1, &  \ldots  &    \G1  - 2,  &
     \G1  - 1,  &     \G1    \\
     \end{array}
\]
 where numbers $\{  2, 4, 6,  \hspace{2mm} \ldots  \hspace{2mm}
\G1-4, \G1-2, \G1 \}$ are even and natural (they are
$\frac{\G1}{2}$) and numbers $ \{ \G1+2, \G1+4,\hspace{2mm}  \ldots
\hspace{2mm} 2\G1-4,  2\G1-2,
  2\G1 \}$ are even and extended natural, they also are
  $\frac{\G1}{2}$. Thus, all the properties (i)--(iii) from Section 2.2 hold and the paradox does not occur.

\subsection{Hilbert's paradox of the Grand Hotel}  Let us consider now the Grand Hotel in the \G1-based framework. In the
paradox, the number of   rooms in the Hotel is infinite. In the new
terminology  it is not sufficient to say this, it is required\emph{ to
indicate explicitly} the infinite number of the rooms in the Hotel. Suppose that it
has \ding{172} rooms. When a new guest arrives, it is proposed to
move the guest occupying room~1 to room~2, the guest occupying room~2 to room~3, etc. At the end of this procedure the guest from the
last room having the number \ding{172} should be moved to the room
\ding{172}+1. However,  the Hotel has only \ding{172} rooms and,
therefore, the poor guy from the room~\G1 will go out of the Hotel
(the situation that would occur  in  hotels with a finite number of rooms,
if such a procedure would be implemented). A nice animation describing
 Hilbert's paradox of the Grand Hotel  in the \G1-based framework can be viewed at the didactic web page \cite{web_teaching} developed at the University of East Anglia, UK.

Notice once again that there is no contradiction between the two ways to see
the Grand Hotel. The traditional answer is that it is possible to put
the newcomer in the first room. The \G1-based way of doing confirms this
result  but shows something that was invisible traditionally -- the
guest from the last room should go out of the Hotel. Thus, the paradox is avoided.

\subsection{Three paradoxes regarding divergent series}

Let us consider the definition of $x$ in (\ref{3.100}). Thanks to \G1, we have different infinite integers and, therefore, we can consider sums having different infinite numbers of summands.
Thus, with respect to the new methodology, (\ref{3.100}) is not well
defined because the number  of summands in the sum (\ref{3.100})
is not explicitly indicated. Recall that to say just that there are $\infty$ many summands has the same meaning of the phrase `There are many summands' (cf.  (\ref{1}),(\ref{2})).

Thus, it is necessary to indicate explicitly an infinite number of addends, $k$, (obviously, it can be finite, as well). After this  (\ref{3.100}) becomes
\[
x(k)=1+2+4+8+\ldots+2^{k-1}
\]
and multiplying both parts by two and adding one to both the right and left sides of the equality gives us
\[
2x(k)+1=1+2+4+8+\ldots+2^{k-1}+2^{k}.
\]
Thus, when we go to substitute, we can see that there remains an addend, $2^{k}$, that is infinite if $k$ is infinite, that was invisible in the traditional framework:
\[
2x(k)+1=\underbrace{1+2+4+8+\ldots+2^{k-1}}_{x(k)}+2^{k}.
\]
The substitution gives us the resulting formula
\[
x(k)=2^{k}-1
\]
that works for both finite and infinite values of $k$ giving different results for different values of $k$ (exactly as it happens for the cases with finite values of $k$). For instance, $x(\G1)=2^{{\tiny\G1}}-1$  and $x(3\G1)=2^{{3\tiny\G1}}-1$. Thus,  the
paradox (\ref{-1}) does not take place.

Let us consider now   Grandi's series. To
calculate the required sum, we should indicate
explicitly the number of addends, $k$, in it. Then it follows
that
 \beq
S(k)=\underbrace{1-1+1-1+1-1+1-\ldots}_{k\,\,\mbox{\small
addends}} = \left \{
\begin{array}{ll} 0, &
  \mbox{if  } k=2n,\\
1, &    \mbox{if  } k=2n+1,\\
 \end{array} \right.
  \label{1-1}
 \eeq
and it is not important whether $k$ is finite or infinite. For
example, $S(\G1)=0$ since  \ding{172}   is even. Analogously, $S(\G1-1)=1$ because $\G1-1$ is
odd.

As it happens in the cases where the  number of addends in a sum is
finite, the result of summation does not depend on the way the
summands  are rearranged. In fact, if we know the exact infinite
number of addends and the order the signs are alternated  is clearly
defined, we know also the exact number of positive and negative
addends in the sum. Let us illustrate this point by supposing, for
instance, that we want to rearrange addends in the sum $S(2\G1)$
as follows
\[
S(2\G1)=1+1-1+1+1-1+1+1-1+\ldots
\]
Traditional mathematical tools used to study divergent series give
an impression that this rearrangement modifies the result. However,
in the \G1-based framework we know that this   is just a consequence
of the weak lens used to observe infinite numbers. In fact, thanks
to \G1 we are able to fix an infinite number of summands. In our
example the sum has $2\G1$ addends, the number $2\G1$ is even and,
therefore, it follows from (\ref{1-1}) that $S(2\G1)=0$. This
means also that in the sum there are $\G1$ positive and $\G1$
negative items. As a result, addition of the groups 1+1--1
considered above can continue only until the positive units present
in the sum will not finish and then there will be necessary to
continue to add only negative summands. More precisely, we have
 \beq
S(2\G1)=\underbrace{1+1-1+1+1-1+\ldots+1+1-1}_{\mbox{\tiny{\G1}}\,\,
 \mbox{\small positive and
}\frac{\mbox{\tiny{\G1}}}{2}\mbox{\small negative
addends}}\hspace{1mm}\underbrace{-1-1-\ldots-1-1}_{\frac{\mbox{\tiny{\G1}}}{2}\mbox{
\small negative addends}}=0,
   \label{1+1-1}
 \eeq
  where the result of the first part in this
rearrangement is calculated as
$(1+1-1)\cdot\frac{\mbox{\tiny{\G1}}}{2}=\frac{\mbox{\tiny{\G1}}}{2}$
and the result of the second part summing up negative units is equal
to $-\frac{\mbox{\tiny{\G1}}}{2}$ giving so the same final result $S(2\G1)=0$. It becomes clear from (\ref{1+1-1}) the origin of the
Riemann series theorem. In fact, the second part of (\ref{1+1-1})
containing only negative units is invisible if one works with the
traditional numeral $\infty$.

Let us use now the \G1 lens to observe Ramanujan's paradoxical result (\ref{Ramanu_3}). As it was in the summation discussed above, it is necessary to indicate explicitly an infinite number of addends, $n$, in the sum
 \beq
 c(n)=1+2+3+ 4+ 5+ \ldots + n.
 \label{Ramanu_7}
 \eeq
The \G1 methodology allows us to compute this sum for infinite
values of $n$ directly  (see   \cite{EMS} for a detailed discussion) and to show that for
infinite (and finite) values of $n$ it follows
 \beq
 c(n)=0.5 n (1+n)
 \label{Ramanu_8}
 \eeq
and by taking $n=\G1$ we can easily compute the sum of all natural numbers
 \beq
  c(\G1)=1+2+3+ 4+ 5+ \ldots + (\G1-2) + (\G1-1) + \G1.
 \label{Ramanu_8.1}
 \eeq
that, obviously, is equal to $0.5 \G1 (1+\G1)$.

Let us now return to Ramanujan summation and consider the main trick
of (\ref{Ramanu_5}) consisting of displacement of addends in its
second line. Since we work with natural numbers, we have $n=\G1$ addends in   the sum (\ref{Ramanu_7}). As a consequence, the
displacement of (\ref{Ramanu_5}) can be re-written more accurately with the observation of the last addends in each line of (\ref{Ramanu_5}) as follows
 \beq
 \begin{array}{rcrrrrrrrrr}
  c(\G1) & = & 1  + 2&\hspace{-3mm} + 3& \hspace{-3mm}+ 4&\hspace{-3mm}+ 5 & \hspace{-3mm} +\ldots &\hspace{-3mm} + \G1-1 &\hspace{-6mm} + \G1,& \\
   4c(\G1) & = &   4& \hspace{-3mm}&\hspace{-3mm} + 8 &\hspace{-3mm} & \hspace{-3mm}+\ldots &\hspace{-3mm} & \hspace{-6mm}+  4 \frac{\G1}{2}& \hspace{-3mm} +4 \left( \frac{\G1}{2}+1\right) + \ldots +
 4 \G1, \\
  -3c(\G1) & = & 1  -2&\hspace{-3mm}+ 3& \hspace{-3mm} -4&\hspace{-3mm}+  5    & \hspace{-3mm}+  \ldots  &\hspace{-3mm}+\G1-1 &
  \hspace{-6mm}-\G1 & \hspace{-3mm} -4 \left( \frac{\G1}{2}+1\right) - \ldots
  - 4 \G1.
   \end{array}
 \label{Ramanu_6}
 \eeq
Thus to   $0.5\G1$ even addends in
(\ref{Ramanu_8.1}) there will be added the first $0.5\G1$ numbers from
(\ref{Ramanu_8.1}) multiplied by 4, i.e., each even number~$i$
from the first line of (\ref{Ramanu_6})  will be summed up with the number
$4\cdot\frac{i}{2}$  whereas $0.5\G1$  odd $i$ from (\ref{Ramanu_8.1}) are summed up with
zeros. The displacement of the second line in  (\ref{Ramanu_6}) leads to the fact that only $0.5\G1$ summands of this line will participate in this addition and there will remain  $0.5\G1$
more addends that were invisible in the traditional framework. They are
 \beq
 4\cdot \left( \frac{\G1}{2}+1\right)+4\cdot\left(\frac{\G1}{2}+2\right)+ \ldots +
 4\cdot\G1.
  \label{Ramanu_9}
 \eeq

Let us compute now the right-hand part of the third line of
(\ref{Ramanu_6})   using  the fact that we can rearrange addends in the sum without changing the result. In such way we can group the addends in three   arithmetical progressions having $0.5\G1$ addends each
 \[
   1  -2 + 3  -4 +  5   + \ldots  +\G1-1  -\G1  -4
\left( \frac{\G1}{2}+1\right) - 4 \left(\frac{\G1}{2}+2\right) -
\ldots
  - 4 \G1 =
 \]
 \[
 \underbrace{(1     + 3   + 5 + \ldots + (\G1-3) + (\G1-1) )}_{=   (1+(\mbox{\tiny{\ding{172}}}-1)) \mbox{\tiny{\ding{172}}}/4  } -  \underbrace{( 2 +  4 + 6 +
\ldots + (\G1-2) + \G1  )}_{=  (2+\mbox{\tiny{\ding{172}}}) \mbox{\tiny{\ding{172}}}/4}-
 \]
 \[ -4\underbrace{\left(
\left( \frac{\G1}{2}+1\right) + \left(\frac{\G1}{2}+2\right)+\ldots
+(\G1-1)+\G1 \right)}_{= ((\frac{\mbox{\tiny{\ding{172}}}}{2}+1)+\mbox{\tiny{\ding{172}}}) \mbox{\tiny{\ding{172}}}/4}=
  \]
 \[
      (1+(\G1-1))\frac{\G1}{4} -  (2+\G1) \frac{\G1}{4}  -4  \left(\left(\frac{\G1}{2}+1\right)
+\G1\right)\frac{\G1}{4}=-3\frac{\G1}{2}\left( \G1 +1\right).
 \]
 Thus, we have obtained that
 \[
  -3c(\G1) = -3\frac{\G1}{2}\left( \G1 +1\right).
  \]
 As was expected, the obtained result shows that  the third line of
(\ref{Ramanu_6}) is consistent with (\ref{Ramanu_8})  for $n=\G1$.

Thus,  it has been shown     that    Riemann's result on
rearrangements of addends in series  is a consequence of the fact
that symbol $\infty$ used traditionally does not allow us to express
quantitatively the infinite number of addends in the series.   The
usage of the grossone methodology allows us to see that (as it
happens in the case where the number of addends is finite)
rearrangements of addends do not change the result for any sum with
a fixed   infinite number of summands. This happens because if one
knows the number of addends and the rule used to alternate their
signs, he/she knows the number of positive and negative addends.
Thus, the careful
counting of the number of addends in infinite series allows us to
avoid this kind of paradoxical results if \G1-based numerals are used.

\subsection{The rectangle paradox of Torricelli} This paradox involves infinitesimals. The numeral system based on \G1 allows us to express them easily and to execute arithmetical operations with them. For instance,  numbers  consisting of addends having    negative
finite  powers of \G1 represent infinitesimals and the   simplest infinitesimal number is
$\mbox{\ding{172}}^{-1}=\frac{1}{\mbox{\ding{172}}}$. It is the inverse element with respect to
multiplication for~\ding{172}:
 \beq
\mbox{\ding{172}}^{-1}\cdot\mbox{\ding{172}}=\mbox{\ding{172}}\cdot\mbox{\ding{172}}^{-1}=1.
 \label{3.15.1}
       \eeq
The following two numbers are other examples of infinitesimals: $5.1\G1^{-2}, -6.1\G1^{-3}+5.1\G1^{-32}$, etc. Note that all infinitesimals are not equal to zero. In particular,
$\frac{1}{\mbox{\ding{172}}}>0$ because it is a result of division
of two positive numbers. It also has a clear granary interpretation. Namely, if we
have a sack   containing \ding{172} seeds, then one sack divided by
the number of seeds in it is equal to one seed. Vice versa, one
seed, i.e., $\frac{1}{\G1}$, multiplied   by the number of seeds in
the sack, $\G1$, gives one sack of seeds.

 \begin{figure}[t]
   \centering
   \includegraphics[width=11.3cm, height=6.5cm]{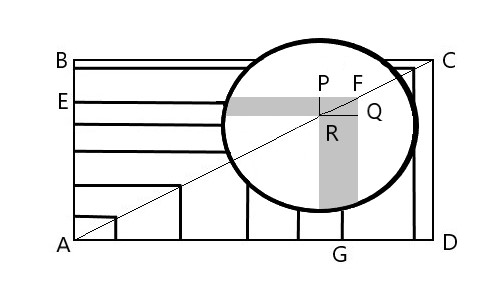}
  \caption{ The rectangle paradox of Torricelli in the grossone-based framework  }
   \label{fig_2}
  \end{figure}

To consider Torricelli's paradox in the grossone framework  it is worthy to mention that sums with an infinite number of infinitesimal addends can give infinitesimal, finite, or infinite results in dependence of the number of summands and their value. As an example, let us consider this sum
 \beq
T(k)=\underbrace{\G1^{-2}+\G1^{-2} +\G1^{-2} \ldots + \G1^{-2}+\G1^{-2}}_{k\,\, \mbox{  addends}}.
   \label{Tk}
 \eeq
Then, for $k= 2\G1$, we obtain an infinitesimal result and for $k=3\G1^2 $  and $k= 4\G1^3 $   finite and infinite results, respectively:
\[
T(2\G1) = \G1^{-2} \cdot  2\G1  = 2\G1^{-1}, \hspace{5mm}   T(3\G1^2) = \G1^{-2} \cdot 3\G1^2 = 3,
 \]
 \[
  T(4\G1^3) = \G1^{-2} \cdot 4\G1^3 = 4\G1.
\]

This machinery allows us to compute directly the areas of the two triangles ABC and CDA  using \G1-based infinitesimals (see Fig.~\ref{fig_2}). In order to be able to execute numerical computations, let us suppose that the length $|AB|=1$, $|BC|=2$ and the line AC has no width (for instance, the rectangle was cut along this line). Obviously, it is also possible to consider the situation where AC has an infinitesimal width but this point is not so important for the essense of the paradox. The procedure that has led to the paradox required to cover the triangles by segments having an infinitesimal width. Without loss of generality and for simplicity suppose that our horizontal segments have the width $h=\G1^{-1}$ (it is easy to see that by taking $h=2\G1^{-1}$ or $h=\G1^{-3}$  the results will be analogous). Then, since $|AB|=1$, it will take \G1 segments of the width $\G1^{-1}$ to cover the whole  triangle ABC. By construction, the triangle CDA will also be covered by \G1 segments. Fig.~\ref{fig_2} considers one horizontal segment, EF, and the corresponding vertical segment, FG, and shows under a magnifying glass the situation in the neighborhood of the point F. The \G1-based framework allows us to observe that both horizontal and vertical segments have triangular ends touching the line AC. Moreover, we can calculate easily the area of these triangles and both the width and the length of the vertical segments.

Since $|AB| /  |BC|= |PR| /  |RQ| $ and $|PR|= \G1^{-1}$, it follows immediately that $|RQ|=2\G1^{-1}$. As a result,   areas of the triangles RPF and FQR (and of other $2\G1-2$ similar triangles on the line AC) are equal to $\G1^{-1}\cdot 2\G1^{-1}/2=\G1^{-2}$. Thus, each horizontal segment $i, 1 \le i \le \G1,$ consists of a rectangle with the width $\G1^{-1}$ and the length $2-2\G1^{-1}i$, having so the area $$\G1^{-1}(2-2\G1^{-1}i)=2\G1^{-1}-2\G1^{-2}i,$$  and of a triangle similar to RPF   having the area equal to $\G1^{-2}$ (notice that for $i=\G1$ the rectangle is absent). Therefore, the area $S_{ABC}^{i}$  of the $i$th horizontal segment is
 \[
 S_{ABC}^{i} = 2\G1^{-1}-2\G1^{-2}i + \G1^{-2}.
 \]
 In order to obtain the area $S_{ABC}$ of the whole coverage of the triangle ABC it is sufficient just to sum up the areas of all \G1 small segments
\[
S_{ABC} = \Sigma_{i=1}^{\G1}  S_{ABC}^{i} = \Sigma_{i=1}^{\G1} (2\G1^{-1}-2\G1^{-2}i + \G1^{-2}) = 2 + \G1^{-1}  - 2\G1^{-2} \Sigma_{i=1}^{\G1}i =
 \]
\[
  2 + \G1^{-1}  - 2\G1^{-2} (\G1+1)\G1/2 = 2 + \G1^{-1}  - (1+\G1^{-1}) = 1.
 \]
The area of the triangle CDA is calculated by a complete analogy.  Each  vertical segment $i, 1 \le i \le \G1,$ consists of a rectangle with the width $2\G1^{-1}$ and the height $1-\G1^{-1}i$, having so the area $$2\G1^{-1}(1-\G1^{-1}i)=2\G1^{-1}-2\G1^{-2}i,$$   and of a triangle similar to FQR   having the area equal to $\G1^{-2}$. Therefore, the area $S_{CDA}^{i}$  of the $i$th vertical segment is
 \[
 S_{CDA}^{i} = 2\G1^{-1}-2\G1^{-2}i + \G1^{-2}= S_{ABC}^{i}.
 \]
Since the number of horizontal and vertical segments is equal, this fact completes the consideration and shows that this paradox is also avoided.

\subsection{Thomson's lamp paradox}

In order to reconsider the Thomson lamp paradox let us
  remind   traditional definitions of
infinite sequences and subsequences.  An \textit{infinite sequence}
$\{a_n\}, a_n \in A, n \in \mathbb{N},$ is a function having as the
domain the set of natural numbers, $\mathbb{N}$, and as the codomain
a set $A$. A \textit{subsequence} is obtained from  a sequence by
deleting  some (or possibly none) of its elements.  In a sequence
$a_1, a_2, \ldots, a_{n-1}, a_n$ the number $n$ is the number of
elements of the sequence. Traditionally, only finite values of $n$
are considered. Grossone-based numerals give us the possibility to
observe infinite numbers and, therefore, to see not only the initial
elements of an infinite sequence $a_1, a_2, \ldots $ but also its
final part  $ \ldots, a_{n-1}, a_n$ where $n$ can assume different
infinite values.

  The IUA  states that the set of natural numbers, $\mathbb{N}$, has
\ding{172} elements. Thus, by the  above definition, any sequence
having $\mathbb{N}$ as the domain  has \ding{172} elements. Since any subsequence is obtained by deleting some (or possibly none) of
the \G1  elements from a sequence, any sequence can have   at most
grossone elements.

Since the switches are executed in a sequence, the maximal number of switches that can be done is equal to \G1. Remind also that we have already established that \G1 is even. Thus, after \G1 switches the lamp will be \emph{off} if   initially it was \emph{on}  and, vice versa, it will be \emph{on} if initially it was \emph{off}.

The \G1-based methodology gives us the opportunity to calculate how much time will take this procedure of switching the lamp. Remind that it is on for $\frac{1}{2}$ minute, then  it is off for $\frac{1}{4}$ minute, then  again on for $\frac{1}{8}$ minute, etc. Thus we deal with the sum of \G1 addends of the form $\frac{1}{2^i}, 1 \le i \le \G1$. It is easy to show (see \cite{EMS} for details) that
\[
\Sigma_{i=1}^{\G1} \frac{1}{2^i} = 1 -\frac{1}{2^{\mbox{\tiny{\G1}}}},
\]
i.e.,   this procedure of switches will not reach number one, it will be infinitesimally  close to one.

\section{A brief conclusion}
 It has been shown in this article that    the surprising  counting systems  of   Amazonian tribes,  Pirah\~{a} and Munduruk\'u, open an interesting perspective on  some classical paradoxes of infinity. The opportunity to use many different numerals to deal with infinities and infinitesimals offered by a recently introduced computational methodology has allowed us to switch from qualitative considerations of paradoxes to their quantitative analysis.

\section*{Declarations}
The author states that there is no conflict of interest and no funding has been used to execute this research.

\end{document}